\documentclass[11pt]{amsart}
\usepackage{amssymb}
\usepackage{amsthm}
\usepackage{amsmath}
\usepackage{amsmath,amssymb,mathrsfs}
\usepackage{xcolor}
\usepackage[T1]{fontenc}
\usepackage[utf8]{inputenc}

\def\e{\varepsilon}

\def\Re{{\rm Re}\,}

\newtheorem{thm}{Theorem}[section]
\newtheorem{lemma}[thm]{Lemma}

\newtheorem{theorem}[thm]{Theorem}

\theoremstyle{definition}

\newtheorem{example}[thm]{Example}

\theoremstyle{remark}

\theoremstyle{definition}

\theoremstyle{remark}

\title[Preserver Problems for Products of Operator Differences]{Preserver Problems for Products of Operator Differences}

\author{Arup Majumdar}
\address{Institute of Mathematics, Czech Academy of Sciences,
\v{Z}itn\'a 25, 115~67 Prague 1, Czech Republic}
\email{arupmajumdar93@gmail.com}

\begin{document}
\begin{abstract}
  In this paper, we investigate surjective maps $\Phi:\mathscr B(H)\to\mathscr B(H)$ preserving various operator quantities associated with products of differences, namely,
$$\alpha\big((\Phi(A)-\Phi(B))\Phi(C)\big)
=\alpha\big((A-B)C\big)
\quad \text{for all}~ A,B,C\in\mathscr B(H),$$
where $\alpha$ denotes one of the following: the numerical radius, maximal numerical range, minimum modulus, surjectivity modulus, reduced minimum modulus, or spectral radius. Under suitable assumptions, we establish structural characterizations of these maps and determine their possible forms.
\end{abstract}

 \maketitle

\noindent
\textbf{Keywords:}
maximal numerical range, surjectivity modulus, spectral radius.

\noindent
\textbf{Mathematics Subject Classification:}
 	15A86, 47A12, 47B49.

    \section{Introduction}

Preserver problems constitute an important and active area of research in operator theory. The main objective of this area is to characterize maps between operator algebras that preserve certain operator quantities, algebraic relations, or spectral properties. In particular, considerable attention has been devoted to preservers of the spectral radius, numerical range, numerical radius, maximal numerical range, minimum modulus, surjectivity modulus, and reduced minimum modulus.

Several authors have investigated preserver problems involving products of operators. For instance, Hou and Di \cite{hou2006maps} studied maps preserving the numerical range of operator products. Zhang, Hou, and He \cite{Xiuling} investigated maps preserving the numerical radius and cross norms of operator products. Related results concerning numerical-radius-preserving maps for products of self-adjoint operators were obtained by He, Hou, and Zhang \cite{he2010maps}. Furthermore, Dhifaoui and Mabrouk \cite{dhifaoui2024maps} characterized maps preserving the maximal numerical range of operator products. For further results concerning the maximal numerical range of triple products of operators, we refer to Bourhim and Mabrouk \cite{bourhim2025maps}.

Preserver problems involving other operator moduli have also received considerable attention. Bourhim, Mashreghi, and Stepanyan \cite{bourhim2014nonlinear} studied nonlinear maps preserving the minimum and surjectivity moduli. In the linear setting, Skhiri \cite{skhiri2009linear} characterized surjective linear maps preserving these moduli on the algebra of bounded operators on an infinite-dimensional complex Hilbert space. Moreover, Hajighasemi and Hejazian \cite{Sepide} investigated surjective maps preserving the reduced minimum modulus of products of operators. Other related contributions include the study of maps preserving the reduced minimum modulus of differences of operators in \cite{JavadMashreghi}, as well as maps preserving the reduced minimum modulus of operator products in \cite{Sepide}.

The aforementioned results provide the main motivation for the present work. In this paper, we investigate a related but different class of preserver problems involving products of operator differences, namely,
$$(A-B)C,\qquad A,B,C\in\mathscr B(H).$$
Such expressions combine the additive and multiplicative structures of the operator algebra, making their preservation problems particularly interesting. In contrast to the preservation of ordinary operator products, the presence of the difference \(A-B\) introduces additional challenges in determining the structure of the preserving maps.

More precisely, we study surjective maps on \(\mathscr B(H)\) satisfying preservation identities involving the numerical radius, reduced minimum modulus, maximal numerical range,  surjectivity modulus, minimum modulus, and spectral radius of products of operator differences. Our aim is to obtain structural characterizations of these maps under suitable assumptions and to determine how the preservation of these operator quantities constrains the form of the underlying transformations. Thus, the present work extends several existing directions in operator preserver theory to products involving differences of operators.\\
The paper is organized as follows. Section 2 is devoted to the basic notation and preliminary results that will be used in the proofs of the main theorems in Section 3. Section 3 investigates the forms of all surjective maps (not necessarily linear, unless otherwise specified) \(\Phi:\mathscr B(H)\to\mathscr B(H)\) satisfying
\[
\alpha\big((\Phi(A)-\Phi(B))\Phi(C)\big)
=\alpha\big((A-B)C\big),
\]
for all \(A,B,C\in\mathscr B(H)\), where \(\alpha\) denotes one of the following quantities: the numerical radius, reduced minimum modulus, maximal numerical radius, surjectivity modulus, minimum modulus, or spectral radius.

\section{Basic Notations}
Throughout this article, \(H\) denotes a finite- or infinite-dimensional complex Hilbert space unless otherwise specified, and $\mathscr{B}(H)$ denotes the algebra of all bounded linear operators on \(H\). The range and null space of $T$ are denoted by $R(T)$ and $N(T)$ respectively. For \(x\in H\) and \(\varepsilon>0\), \(B(x,\varepsilon)\) denotes the open ball in \(H\) centered at \(x\) with radius \(\varepsilon\). For any subset \(W\) of \(H\), we denote by \(\overline{W}\) its closure in \(H\). $\mathbb{T}$ is the unit cycle of the complex plane. \(A^*\) denotes the adjoint of an operator \(A\in\mathscr{B}(H)\). For any \(x,y\in H\), the rank-one operator \(x\otimes y\) on \(H\) is defined by
\[
(x\otimes y)z=\langle z,y\rangle x,
\qquad \text{for all}~ z\in H.
\]

The numerical range and numerical radius of \(A\) are defined,
respectively, by
\[
W(A)=\{\langle Ax,x\rangle:
x\in H,\ \|x\|=1\},
\]
and
\[
w(A)=\sup\{|\lambda|:\lambda\in W(A)\}.
\]

The maximal numerical range and the maximal numerical radius of
\(A\) are defined, respectively, by
\[
W_0(A)=  \{\lambda \in \mathbb{C}: \langle A x_{n}, x_{n}\rangle \to \lambda \text{ where } \|x_{n}\|= 1 \text{ and } \|A x_{n}\| \to \|A\|\},
\]
and
\[
w_0(A)=\sup\{|\lambda|:\lambda\in W_0(A)\}.
\]

The minimum modulus  and the surjectivity modulus of \(A\) are defined, respectively, by
$$m(A)=\inf\{\|Ax\|:
x\in H,\ \|x\|=1\}~ \text{and}~
q(A)=\sup\{\varepsilon\geq 0:
\varepsilon B(0, 1)\subseteq A(B(0,1))\}.$$ The spectral radius of an operator \(A\in\mathscr B(H)\) is defined by
\[
r(A)=\sup\{|\lambda|:\lambda\in\sigma(A)\},
\]

where \(\sigma(A)\) denotes the spectrum of \(A\). Finally, the reduced minimum modulus of \(A\) is defined by
\[
\gamma(A)=\inf\{\|Ax\|:
x\in N(A)^\perp,\ \|x\|=1\},
\]
with the convention \(\gamma(0)=+\infty\). Moreover, \(A\sim_u B\) means that \(A\) and \(B\) are unitarily similar; that is, there exists a unitary operator \(U\in\mathscr{B}(H)\) such that \(A=UBU^*\). Let \(\{e_\alpha\}_{\alpha\in\Lambda}\) be a fixed but arbitrary orthonormal
basis of \(H\). The conjugation \(J:H\to H\) associated with
this basis is defined by
\[
J\left(\sum_{\alpha\in\Lambda}x_\alpha e_\alpha\right)
=
\sum_{\alpha\in\Lambda}\overline{x_\alpha}\,e_\alpha.
\]
The transpose of \(A\in B(H)\), with respect to this fixed
orthonormal basis, is defined by $A^{\mathrm{tr}}=JA^*J.$

\section{Main Theorems}
\begin{theorem}
Let $\Phi: \mathscr{B}(H) \to \mathscr{B}(H)$ be a surjective map. Then 
$$\|(\Phi(A)-\Phi(C))B\| = \|(A - C)B\|, \text{ for all } A,B, C \in \mathscr{B}(H)$$ if and only if there exists a unitary operator $U \in \mathscr{B}(H)$ and fixed $T_0 \in \mathscr{B}(H)$ such that $$\Phi(A) = UA +T_0, \text{ for all } A \in \mathscr{B}(H).$$
\end{theorem}
\begin{proof}
Taking $B= I,$ we obtain $$\|\Phi(A) - \Phi(C)\| = \|A-C\|,$$ for all $A, C \in \mathscr{B}(H).$ By The Mazur-Ulam Theorem, there exists a surjective real linear isometry operator $F$ on $\mathscr{B}(H)$ such that $$\Phi(X) = F(X) + \Phi(0),$$ for all $X \in \mathscr{B}(H).$ Now taking $C= 0,$ we get $\|F(A)B\| = \|AB\|$ for all $A, B \in \mathscr{B}(H).$ For two non-zero vectors $x, y \in H,$ the one rank operator $B = x \otimes y$ gives $\|F(A)x \otimes y\| = \|Ax \otimes y\|.$ So $\|F(A)x\| = \|Ax\|$ for all $x \in H.$

For each non zero $x \in H,$ Define $U_x : H \to H$ by $$U_x(Ax) = F(A)x,$$ for all $A \in \mathscr{B}(H).$\\
For $A,A_0 \in \mathscr{B}(H)$ with $Ax = A_0 x$ implies $\|F(A-A_0)x\| = 0.$ Thus $U_x$ is well-defined. Moreover, $\|U_x(Ax)\| = \|Ax\|$ gives $U_x$ is a real linear isometry. Now we show that $U_x$ is surjective. Since $F$ is surjective. Then given $w \in H,$ there exists a linear operator $D \in \mathscr{B}(H)$ and $F(E) =D$ for some $E \in B(H),$ such that $$w = Dx = F(E)x = U_x(Ex).$$ 
It follows that $U_x$ is surjective real linear isometry.
Now we divide two cases to prove the statement.\\
\textit{Case 1:} Consider $\dim(H) >1.$ Let $x, y \in H$ be two linearly independent vectors. For arbitrary $u \in H,$ choose $A \in \mathscr{B}(H)$ satisfying $$Ax = u \text{ and } Ay =0.$$ Then 
\begin{align*}
U_{x +y}(u) &= U_{x+y}(A(x+y))\\ 
&= F(A)(x+y)\\
&= U_x(Ax) + U_y(Ay)\\
&= U_x(u).
\end{align*}
Hence $U_{x+y} = U_x.$ Similarly, We can show that $U_{x+y} =U_y.$ Consequently, $U_x =U_y.$ Write their common value as $U,$ we have $$F(A)x = U(Ax),$$ for all $A \in \mathscr{B}(H),$ where $U$ is real linear surjective isometry. Thus $F(A) = UA,$ for all $A \in \mathscr{B}(H).$ Taking $A = I,$ we obtain $U = F(I) \in \mathscr{B}(H)$ is surjective linear(complex) isometry. So $U$ is unitary. Therefore, $$\Phi(A) = UA + \Phi(0).$$

\textit{Case 2:} Consider $\dim(H) =1.$  Then we can directly define $U$ instead of $U_x$ to get the identity $$\Phi(A) = UA + \Phi(0),$$ where $U$ is unitary.\\

Conversely, Suppose $\Phi(A) = UA +T_0,$ where $U$ is unitary and $T_0 \in \mathscr{B}(H).$ Then it is evident that $$\|(\Phi(A) - \phi(C))B\| = \|(A - C)B\|,$$ for all $A, B, C \in \mathscr{B}(H).$
\end{proof}
A natural question is whether the numerical radius distance alone determines the algebraic and geometric structure of such maps, without assuming linearity or continuity. Bai and Hou answered this question by characterizing all surjective numerical radius distance-preserving maps on $\mathscr{B}(H)$. Their result shows that every such map is built from unitary equivalence, transpose, adjoint, multiplication by a scalar of modulus one, and translation by a fixed operator.

\begin{theorem}\cite{bai2004numerical}\label{thm 3.2}
Let $\Phi: \mathscr{B}(H)\to \mathscr{B}(H)$ be a surjective map. Then
$$w(\Phi(A)-\Phi(B))=w(A-B),
\quad \text{for all}~ A,B\in \mathscr{B}(H),$$
if and only if there exists a unitary operator $U\in \mathscr{B}(H)$,
a complex scalar $\mu$ with $|\mu|=1$, and an operator
$S\in \mathscr{B}(H)$ such that one of the following holds:
$$\begin{aligned}
\text{\rm (i)}\quad
&\Phi(A)=\mu UAU^*+S \text{ for all } A\in \mathscr{B}(H),\\
\text{\rm (ii)}\quad
&\Phi(A)=\mu UA^{\mathrm{tr}}U^*+S \text{ for all } A\in \mathscr{B}(H),\\
\text{\rm (iii)}\quad
&\Phi(A)=\mu U(A^*)^{\mathrm{tr}}U^*+S \text{ for all } A\in \mathscr{B}(H),\\
\text{\rm (iv)}\quad
&\Phi(A)=\mu UA^*U^*+S \text{ for all } A\in \mathscr{B}(H).
\end{aligned}$$
Here \(A^{\mathrm{tr}}\) denotes the transpose of \(A\)
with respect to a fixed but arbitrary orthonormal basis
of $H.$
\end{theorem}
Theorem \(\ref{thm 3.3}\) establishes the structural characterization of surjective mapping $\Phi: \mathscr{B}(H) \to \mathscr{B}(H)$ that preserve the numerical range of sums of operators.

\begin{theorem}\label{thm 3.3}
Let $\Phi: \mathscr{B}(H)\to \mathscr{B}(H)$ be a surjective map. Then $$W(\Phi(A) + \Phi(B)) = W(A +B)$$ for all $A, B \in \mathscr{B}(H)$ if and only if  if there exists a unitary operator $U\in \mathscr{B}(H)$ such that one of the following holds:
$$\begin{aligned}
\text{\rm (i)}\quad
&\Phi(A)= UAU^* \text{ for all } A\in \mathscr{B}(H),\\
\text{\rm (ii)}\quad
&\Phi(A)= UA^{\mathrm{tr}}U^* \text{ for all } A\in \mathscr{B}(H),
\end{aligned}$$
\end{theorem}
\begin{proof}
Taking $A = B = 0.$ Then $W(\Phi(0)) = 0$ which implies that $\Phi(0) = 0.$ It is observed that  $\Phi(-A) = -\Phi(A)$ because $W(\Phi(A) + \Phi(-A)) = 0.$ Taking the numerical radius both sides, we obtain 
$$w(\Phi(A) - \Phi(B)) = w(A- B) \qquad \text{ for all } A,B \in \mathscr{B}(H).$$ Putting $A= 0,$ from Theorem \ref{thm 3.2}, we get $S= 0.$ Consequently we have one of the following conditions:
$$\begin{aligned}
\text{\rm (i)}\quad
&\Phi(A)=\mu UAU^* \text{ for all } A\in \mathscr{B}(H),\\
\text{\rm (ii)}\quad
&\Phi(A)=\mu UA^{\mathrm{tr}}U^* \text{ for all } A\in \mathscr{B}(H),\\
\text{\rm (iii)}\quad
&\Phi(A)=\mu U(A^*)^{\mathrm{tr}}U^* \text{ for all } A\in \mathscr{B}(H),\\
\text{\rm (iv)}\quad
&\Phi(A)=\mu UA^*U^* \text{ for all } A\in \mathscr{B}(H).
\end{aligned}$$
Now putting $B = 0$ in the original condition, we obtain $W (\Phi(A)) = W(A)$ for all $A \in \mathscr{B}(H).$ Now we take $A = I.$ Then $\{\mu\} = W(\Phi(I)) = W(I) = \{1\}.$ So, the four possibilities reduce to 
$$\begin{aligned}
\text{\rm (i)}\quad
&\Phi(A)= UAU^* \text{ for all } A\in \mathscr{B}(H),\\
\text{\rm (ii)}\quad
&\Phi(A)= UA^{\mathrm{tr}}U^* \text{ for all } A\in \mathscr{B}(H),\\
\text{\rm (iii)}\quad
&\Phi(A)= U(A^*)^{\mathrm{tr}}U^* \text{ for all } A\in \mathscr{B}(H),\\
\text{\rm (iv)}\quad
&\Phi(A)= UA^*U^* \text{ for all } A\in \mathscr{B}(H).
\end{aligned}$$
Choose $A = iI.$ Then the four possibilities are $$\text{\rm (i)} ~\Phi(iI) = iI, \text{\rm (ii)}~ \Phi(iI) = iI, \text{\rm (iii)}~ \Phi(iI) = -iI,  \text{\rm (iv)}~ \Phi(iI) = -iI.$$
Consequently, the condition $\text{\rm (iii)}$ and $\text{\rm (iv)}$ are impossible because $W(\Phi (A)) = W(A)$ for all $A\in \mathscr{B}(H).$ Therefore, we left with 
$$\begin{aligned}
\text{\rm (i)}\quad
&\Phi(A)= UAU^* \text{ for all } A\in \mathscr{B}(H),\\
\text{\rm (ii)}\quad
&\Phi(A)= UA^{\mathrm{tr}}U^* \text{ for all } A\in \mathscr{B}(H),
\end{aligned}$$
where $U$ is unitary in $\mathscr{B}(H).$\\
Conversely, suppose one the following holds
$$\begin{aligned}
\text{\rm (i)}\quad
&\Phi(A)= UAU^* \text{ for all } A\in \mathscr{B}(H),\\
\text{\rm (ii)}\quad
&\Phi(A)= UA^{\mathrm{tr}}U^* \text{ for all } A\in \mathscr{B}(H),
\end{aligned}$$
where $U$ is unitary in $\mathscr{B}(H).$ Then it is evident that $W(\Phi(A) + \Phi(B)) = W(A +B),$ for all $A, B \in \mathscr{B}(H).$
\end{proof}
 Hou and He \cite{KhejHou} studied surjective nonlinear maps on the space of self-adjoint operators on a complex separable Hilbert space, characterizing those that preserve the numerical radius and the numerical range of Lie products. Their results provide structural classifications without assuming linearity or continuity of the maps. Theorem \ref{thm 3.4} characterizes numerical radius preserving maps for a different operator products other than Lie products.
 \begin{theorem}\label{thm 3.4}
 Let $\dim H > 2$ and $\Phi: \mathscr{B}(H) \to \mathscr{B}(H)$ be a surjective map. Then
 $$w((\Phi(A)-\Phi(B))\Phi(C)) = w((A-B)C) \text{ for all } A,B,C \in \mathscr{B}(H)$$ if and only if there exists a unitary or conjugate unitary operator $U$ and a complex scalar $\mu$ with $\vert \mu \vert = 1$  such that $$\Phi(A) = \mu UAU^* \text{ for all } A \in \mathscr{B}(H).$$
 \end{theorem}
\begin{proof}
Since $\Phi$ is surjective. Then there exists an operator $C_0 \in \mathscr{B}(H)$ such that $\Phi(C_0) = 0.$ Putting $C= C_0$ in the original equation, we obtain $w((A-B)C_0) = 0,$ for all $A, B \in \mathscr{B}(H).$ So $w(XC_0) = 0$ for all $X \in \mathscr{B}(H).$\\
Suppose $C_0 \ne 0.$ Choose $y \in H$ such that ${C_0}^*y \ne 0$ and put $x = {C_0}^*y$ then take $X = x \otimes y,$ we obtain $$XC_0 = {C_0}^* y \otimes {C_0}^* y = x \otimes x.$$ But $w(XC_0) = \|x\|^2$ is nonzero which is a contradiction. Hence $C_0 = 0.$ Consequently, $\Phi(0) =0.$ Now we show that $\Phi$ is injective. Consider $\Phi(A) = \Phi(B)$ for some $A, B \in \mathscr{B}(H).$ From the original equation, we get 
$w((A-B)C)= 0$ for all $C \in \mathscr{B}(H).$ If $X= A-B \ne 0,$ then choose $u \in H$ and $Xu \ne 0.$ Take $C = u \otimes Xu,$ we have $XC = Xu \otimes Xu.$ So $w((A-B)C) = \|Xu\|^2$ is nonzero which is a contradiction. Thus $A-B =0,$ hence $\Phi $ is bijective.\\
Now we claim that $\Phi$ is real linear. Since $\Phi$ is surjective. There exists an operator $D \in \mathscr{B}(H)$ such that $\Phi(D) = I.$ Taking $C= D$ in the original equation, we have $w(\Phi(A)- \Phi(B)) = w((A-B)D).$

Define a new norm on $\mathscr{B}(H)$ by $$\|X\|_D = w(XD) \text{ for all } X \in \mathscr{B}(H).$$ Indeed, homogeneity and the triangular inequality follow from the corresponding properties of numerical radius. If $\|X\|_D = 0,$ then in the property $w(\Phi(A)- \Phi(B)) = w((A-B)D)$ we take $A = X, B= 0,$ we get $w(\Phi(X)) = \|X\|_D =0.$ It follows from the injectivity of $\Phi$ that $X=0.$ There is another norm space $(\mathscr{B}(H), \|.\|_w)$ defined by $\|T\|_w = w(T)$ for all $T \in \mathscr{B}(H).$ Thus, $\Phi$ is surjective isometry from the real norm space $(\mathscr{B}(H), \|.\|_D)$ to the real norm space $(\mathscr{B}(H), \|.\|_w).$
 By the Mazur-Ulam Theorem $\Phi$ is affine over $\mathbb{R}.$ Since $\Phi(0) = 0.$ It follows that $\Phi$ is real linear.\\
 Putting $B =0$ in the original equation and using $\Phi(0) =0$ gives $$w(\Phi(A)\Phi(C)) = w(AC),$$ for all $A, C \in \mathscr{B}(H).$ By using Theorem 2.1 \cite{Xiuling}, there exists a functional $h :\mathscr{B}(H) \to \mathbb{T}$ and a unitary or conjugate unitary operator $U$ on $H$ such that $$\Phi(T) = h(T)UTU^*$$ for all $T \in \mathscr{B}(H).$ Since $\Phi$ is additive. Then  for all $A, B \in \mathscr{B}(H)$
 \begin{align*}
 \Phi(A + B) &= \Phi(A) + \Phi(B)\\
 (h(A+B)-h(A))A &= (h(A+B)-h(B))B.
 \end{align*}
 When $A, B$ are linearly independent, then $h(A) =h(A+B) =h(B).$ Moreover, when $A, B$ are linearly dependent, then choose a linear operator $E \in \mathscr{B}(H)$ which is linearly independent of both $A, B$ because of $\dim H >2.$ Consequently, $$h(A) = h(A+E) =h(E) = h(E+B) = h(B).$$
 Therefore, $h$ is constant. It shows that $h(T) = \mu \in \mathbb{T}$ for all $T \in \mathscr{B}(H).$ We conclude that $$\Phi(T) = \mu UTU^*,$$ for all $T \in \mathscr{B}(H)$ and $U$ is a unitary or conjugate unitary operator with $\vert \mu \vert = 1.$\\
 The reverse inclusion is obvious because of $w(T) = w(UTU^*),$ for all $T \in \mathscr{B}(H)$ and $U$ is unitary or conjugate unitary on $H.$
\end{proof}

In \cite{JavadMashreghi}, the authors established a structural characterization of bicontinuous bijective maps \(\Phi:\mathscr{B}(\mathscr{X})\to\mathscr{B}(\mathscr{Y})\) preserving the reduced minimum modulus of differences, that is,
\[
\gamma(\Phi(A)-\Phi(B))=\gamma(A-B),
\]

where \(\mathscr{X}\) and \(\mathscr{Y}\) are infinite-dimensional complex Banach spaces. Furthermore, in \cite{Sepide}, the authors characterized surjective maps \(\Phi:\mathscr{B}(H)\to\mathscr{B}(H)\) preserving the reduced minimum modulus of operator products, namely,
\[
\gamma(\Phi(A)\Phi(B))=\gamma(AB),
\qquad A,B\in\mathscr{B}(H),
\] when $\dim H >2.$\\
Theorem \ref{thm 3.5} provides a structural characterization of surjective maps \(\Phi:\mathscr{B}(H)\to\mathscr{B}(H)\) that preserve the reduced minimum modulus of operator products involving differences, namely,
\[
\gamma\big((\Phi(A)-\Phi(B))\Phi(C)\big)
=
\gamma\big((A-B)C\big),
\]

for all \(A,B,C\in\mathscr{B}(H)\).
\begin{theorem}\label{thm 3.5}
Let $\Phi: \mathscr{B}(H) \to \mathscr{B}(H)$ be a surjective map. If $$\gamma\big((\Phi(A)-\Phi(B))\Phi(C)\big)
=
\gamma\big((A-B)C\big)$$ for all $A,B,C\in\mathscr{B}(H),$ then $\gamma(\Phi(A)\Phi(B)) = \gamma(AB)$ for all $A, B \in \mathscr{B}(H).$
\end{theorem}
\begin{proof}
Since $\Phi$ is surjective. There exists a unique $C_0 \in \mathscr{B}(H)$ such that $\Phi(C_0) =0.$ Put $C= C_0$ in the original equation, we obtain $\gamma((A-B)C_0) = 0$ for all $A, B \in \mathscr{B}(H).$ Since $A, B \in \mathscr{B}(H)$ are arbitrary. It follows that $\gamma(XC_0) =0$ for all $X \in \mathscr{B}(H).$ Now we show that $C_0 = 0.$\\
Suppose $C_0 \ne 0.$ Choose $y \in H$ such that $C_0^* y \ne 0$ and take any non zero element $x \in H$ with a rank one operator $X = x \otimes y.$ Then $XC_0 = x \otimes C_0^* y.$ Recall that for a non zero rank one operator $\gamma(u \otimes v) = \|u\|\|v\|.$ Consequently, $\gamma(XC_0) = \|x\|\|C_0^*y\|\ne 0$ is a contradiction. Hence $C_0 =0$ which shows that $\Phi(0) =0.$ Put $B= 0,$ we get $$\gamma(\Phi(A)\Phi(C)) = \gamma(AC)$$ for all $A, C \in \mathscr{B}(H).$
\end{proof}
The converse of Theorem \ref{thm 3.5} is not valid in general.
\begin{example}
Define $\Phi$ on $\mathscr{B}(H)$ by $\Phi(T) = e^{i\|T\|}T$ for all $T \in \mathscr{B}(H).$ Then its inverse is $\Phi^{-1}(T) = e^{-i\|T\|}T.$ Moreover, $\gamma(\Phi(T)\Phi(S)) = \gamma(TS)$ for all $T, S \in \mathscr{B}(H).$ Now consider $A = I, B= 2I, \text{ and } C= I.$ We obtain $$\gamma((\Phi(I)-\Phi(2I))\Phi(I)) = \gamma((e^i - 2e^{2i})I) = \vert 1- 2e^i\vert > 1= \gamma(-I)= \gamma((I-2I)I).$$
\end{example}
In \cite{dhifaoui2024maps}, the authors established a structural characterization of surjective maps $\Phi : \mathscr{B}(H) \to \mathscr{B}(H)$ preserving the maximal numerical range of operator products, namely,
$$W_0(\Phi(A)\Phi(B)) = W_0(AB),$$ for all $A, B \in \mathscr{B}(H)$ if and only if there exists a unitary operator $U \in \mathscr{B}(H)$ such that $\Phi(A) = \pm UAU^*$ for all $A \in \mathscr{B}(H).$ For results concerning preservers of the maximal numerical range of triple operator products, we refer the reader to \cite{bourhim2025maps}.\\
Theorem \ref{thm 3.6} provides a structural characterization of surjective maps $\Phi: \mathscr{B}(H) \to \mathscr{B}(H)$ that preserve the maximal numerical range of operator products involving differences, namely, 
$$W_0((\Phi(A) - \Phi(B))\Phi(C)) = W_0((A-B)C),$$ for all $A, B, C \in \mathscr{B}(H).$
\begin{theorem}\label{thm 3.6}
Let $\Phi : \mathscr{B}(H) \to \mathscr{B}(H)$ be a surjective map. Then $$W_0((\Phi(A) - \Phi(B))\Phi(C)) = W_0((A-B)C),$$ for all $A, B, C \in \mathscr{B}(H)$ if and only if there exists a unitary operator $U$ such that $\Phi(A) = \pm UAU^*$ for all $A\in \mathscr{B}(H).$
\end{theorem}
\begin{proof}
Since $\Phi$ is surjective. There exists $C_0 \in \mathscr{B}(H)$ with $\Phi(C_0) = 0.$ Putting $C =C_0$ in the original equation, we obtain $\{0\}= W_0(0) = W_0((A-B)C_0)$ for all $A, B \in \mathscr{B}(H).$ So $W_0(XC_0) =\{0\}$ for all $X \in \mathscr{B}(H).$ Now we show that $C_0 = 0.$ Suppose $C_0 \ne 0.$ Choose $y \in H$ such that $C_0^* y =x \ne 0$ for some $x \in H.$ Take $X = x \otimes y,$ we have $XC_0 = x \otimes C_0^*y = x \otimes x.$ It follows that $W_0(XC_0) = \{\|x\|^2\} = \{0\}$ which is a contradiction. Consequently, $C_0 = 0.$ Now put $B= 0$ in the original equation, we get $$W_{0}(\Phi(A)\Phi(C)) = W_0(AC),$$ for all $A, C \in \mathscr{B}(H).$ Using Theorem 1.1 \cite{dhifaoui2024maps}, the exists a unitary operator $U \in \mathscr{B}(H)$ such that $\Phi(A) = \pm UAU^*$ for all $A \in \mathscr{B}(H).$\\

Conversely, suppose $\Phi(T) = \pm UTU^*$ for some unitary $U \in \mathscr{B}(H)$ and all $T \in \mathscr{B}(H).$ By using the property $W_0(T) = W_0(UTU^*)$ for all $T \in \mathscr{B}(H),$ we obtain  $$W_0(\Phi(A) - \Phi(B)\Phi(C)) = W_0((A-B)C),$$ for all $A, B, C \in \mathscr{B}(H).$
\end{proof}
In paper \cite{bourhim2014nonlinear}, the authors established a structural characterization of surjective maps $\Phi : \mathscr{B}(\mathscr{X}) \to \mathscr{B}(\mathscr{Y}),$ where $\mathscr{X}, \mathscr{Y}$ both are infinite dimensional complex Banach spaces, preserving the surjectivity (minimum) modulus of operator additions and subtractions, namely,
$$q(\Phi(A) \pm \Phi(B)) = q(A\pm B) (\text{ or } m(\Phi(A) \pm \Phi(B)) = m(A \pm B)),$$ for all $A, B \in \mathscr{B}(\mathscr{X}).$ Skhiri \cite{skhiri2009linear} proved that, for a surjective linear map \(\phi:\mathscr B(H)\to\mathscr B(H)\), preservation of the minimum modulus is equivalent to preservation of the surjectivity modulus, and both characterize maps of the form
\[
\phi(T)=UTV,
\]

where \(U,V\in\mathscr B(H)\) are unitary operators and for all $T \in \mathscr B(H).$\\
Theorem \ref{thm 3.8} provides a structural characterization of surjective  liner maps $\Phi: \mathscr B(H) \to \mathscr{B}(H),$ where $H$ is infinite-dimensional, that preserve the surjectivity modulus of operator products involving differences, namely,
$$q((\Phi((A)-\Phi(B))\Phi(C)) = q((A-B)C),$$ for all $A, B \in \mathscr B(H).$ We first establish a lemma needed for the proof of Theorem \(\ref{thm 3.8}\).
\begin{lemma}\label{lemma 3.8}
Let $A, B \in \mathscr{B}(H)$ be invertible satisfying $$q(AXB) = q(X),$$ for all $X \in \mathscr{B}(H).$ Then there exist unitary operators $U_1, V_1$ and non zero scalars $a,b$ such that $$A= a V_1, \quad B= bU_1, \text{ and} \quad  ab=1.$$
\end{lemma}
\begin{proof}
For invertible operator $X,$ we have $q(X) = \frac{1}{\|X^{-1}\|}.$ Consequently, the original equation gives $$\|B^{-1}YA^{-1}\| = \|Y\|,$$ for all invertible $Y \in \mathscr{B}(H).$ We can approximate rank one operator in $\mathscr{B}(H)$ by invertible rank one operator. Suppose $x \otimes y$ is a rank one operator where non zero $x , y \in H.$ For sufficient small $|\langle x, y\rangle| > \e >0,$ we have an invertible operator $(x \otimes y)+ \e I.$ Choose $|\langle x,y\rangle|> \e_n >0$ and $\e_n \to 0.$ So the sequence of invertible operators $Y_n = (x \otimes y)+ \e_n \to x \otimes y,$ as $n \to \infty.$ Consequently, we get $$ \|B^{-1}(x \otimes y)A^{-1}\|=\lim_{n \to \infty}\|B^{-1}Y_nA^{-1}\|= \lim_{n \to \infty}\|Y_n\|= \|x\otimes y\| = \|x\| \|y\|.$$ It follows that $\|B^{-1}x\|\|(A^{-1})^*y\| = \|x\| \|y\|$ for all non zero $x, y \in H.$\\
So, there is a constant, say $c>0,$ such that $$\|B^{-1}x\| = c\|x\| \quad \text{ and } \quad  \|(A^{-1})^* x\| = \frac{1}{c}\|x\|,$$ for all $x \in H.$ Now, define $U = \frac{1}{c}B^{-1}$ is surjective and isometry, hence $U$ is unitary. Renaming $c^{-1} =b $ and $U^* = U_1,$ we get $B = b U_1.$ Again, define $V = c(A^{-1})^*$ is surjective and isometry, hence $V$ is unitary. Renaming $c =a $ and $V= V_1,$ we obtain $A = aV_1$ with $ab=1.$
\end{proof}
\begin{theorem}\label{thm 3.8}
Let $H$ be an infinite-dimensional complex Hilbert space and $\Phi: \mathscr{B}(H) \to \mathscr{B}(H)$ be surjective linear map.  Then $$q((\Phi(A)- \Phi(B))\Phi(C)) = q((A-B)C),$$ for all $A, B, C \in \mathscr{B}(H)$ if and only if there exists a unitary $U$ and a complex scalar $\lambda$ with $\vert \lambda \vert=1$ such that $\Phi(A) = \lambda UAU^*$ for all $A \in \mathscr{B}(H).$
\end{theorem}
\begin{proof}
Taking $B = -A$ in the main equation and using $q(2T) = 2q(T),$ we get $q(\Phi(A)\Phi(C)) = q(AC)$ for all $A, C \in \mathscr{B}(H).$ Recall that $q(T) >0$ if and only if $T$ is surjective, for all $T \in \mathscr{B}(H).$ Now, we show that $A$ is surjective if and only if so $\Phi(A)$ is.\\
When $A$ is surjective, choose $C_0$ such that $AC_0=I.$ Then $$q(\Phi(A)\Phi(C_0)) = q(AC_0)= q(I)=1.$$ Consequently, $\Phi(A)\Phi(C_0)$ is surjective. So, $\Phi(A)$ is surjective. Similarly, $\Phi(A)$ is surjective implies $A$ is surjective. Since $\Phi$ is linear. Then $A-C$ is surjective if and only if $\Phi(A) - \Phi(C)$ is surjective for all $A, C \in \mathscr{B}(H).$ By using Theorem 4.2 \cite{bourhim2014nonlinear}, there are bounded invertible maps $T, S$ on $H$ either both linear or both conjugate linear such that 
\begin{align*}
\Phi(A) &= TAS + \Phi(0), \text{ for all } A \in \mathscr{B}(H)\\
\Phi(A) &= TA^*S +\Phi(0),  \text{ for all } A \in \mathscr{B}(H),
\end{align*}
where $\Phi(0) =0.$ Now, we eliminate the possibilities involving the alternatives involving conjugate linear factors $T$ and $S.$\\
For all $A \in \mathscr{B}(H),$ 
$$\Phi(A) = TJ(A^*)^{\mathrm{tr}}JS = P(A^*)^{\mathrm{tr}}Q \text{ or } \Phi(A) = TJ A^{\mathrm{tr}} JS = PA^{\mathrm{tr}}Q,$$ where $P= TJ $ and $Q = JS$ both are bounded invertible linear operators in $\mathscr{B}(H).$ Let $V$ be unilateral shift on $H,$ extended by the identity on the orthogonal complement if necessary. Then $V$ is not surjective but $V^{\mathrm{tr}} = JV^*J$ is surjective, since $V^*$ is surjective. It shows that $\Phi(V)$ is surjective when $V$ is not surjective which is a contradiction. Similarly, we eliminate the possibility $\Phi(A) = TA^*S$ for all $A \in \mathscr{B}(H).$  Therefore, the only possible representation is $\Phi(A) = TAS$ for all $A \in \mathscr{B}(H)$ and $T, S \in \mathscr{B}(H)$ are invertible linear operators.\\
So $q(\Phi(A)\Phi(C)) = q(TASTCS) = q(AC)$ for all $A, C \in \mathscr{B}(H).$ Set $C=I,$ we get $q(TASTS) = q(A)$ for all $A \in \mathscr{B}(H).$  Using Lemma \ref{lemma 3.8}, we obtain
$$T = aU, \quad STS = b V, \quad \text{and}\quad ab=1,$$ where $U$ and $V$ are unitaries. Next set $A=I,$ we obtain $q(TSTCS) =q(C)$ for all $C \in \mathscr{B}(H).$ Hence $$TST =a_1U_1, \quad S = b_1 V_1,\quad \text{and} \quad a_1b_1 =1,$$ where $U_1, V_1 \in \mathscr{B}(H)$ are unitaries. Consequently, $T = aU, \quad  S= b_1V_1$ and $bV = ab_1^{2}V_1 UV_1.$ Taking norm in both sides, we get $a^2b_1^{2} = ab =1.$ So $|ab_1| =1.$\\
It follows that $\Phi(A) = ab_1UAV_1,$ where $U$ and $V_1$ both are unitaries. Again 
\begin{align*}
q(a^2b_1^2UAV_1 UCV_1) &= q(AC)\\
q(ARC) = q(AC),
\end{align*}
where $R = V_1U$ is unitary. Take inverse of $A, C$ and use the property $q(X) = \frac{1}{\|X^{-1}}\|$ for invertible operator $X \in \mathscr{B}(H),$ yields $\|C^{-1}R^{-1}A^{-1}\|= \|C^{-1}A^{-1}\|.$ Renaming $Y= C^{-1}$ and $Z= A^{-1},$ we have 
\begin{align}\label{equ 1}
\|YR^{-1}Z\| = \|YZ\|
\end{align}for all invertible operators $Y, Z \in \mathscr{B}(H).$ Now, fixed $Y$ and put $W =YZ,$ the identity (\ref{equ 1}) becomes $\|YR^{-1}Y^{-1}W\| = \|W\|,$ for all invertible $W \in \mathscr{B}(H).$ By the same rank one approximation argument used in Lemma \ref{lemma 3.8} implies that $YR^{-1}Y^{-1}$ is isometry. Since it is invertible. Thus $YR^{-1}Y^{-1}$ is unitary.\\
For every invertible $Y \in \mathscr{B}(H),$ we have 
\begin{align*}
(YR^{-1}Y^{-1})^*(YR^{-1}Y^{-1}) &= I\\
R(Y^*Y) = (Y^*Y)R.
\end{align*}
Every positive invertible operator is of the form $Y^*Y.$ Hence $R$ commutes with every positive invertible operator in $\mathscr{B}(H)$. It follows that $R$ commutes with every self-adjoint operator in $\mathscr{B}(H).$ Since the center of $\mathscr{B}(H)$ consists solely the scalar multiples of the identity. Consequently, $R = \lambda I,$ for some $\lambda \in \mathbb{C}$ with $|\lambda| =1$ because $R$ is unitary and $V_1 = \lambda U^*.$ Therefore, $\Phi(A) =\lambda ab_1UAU^{*}= \mu UAU^*,$ where $|\mu| = |ab_1| =1$ and for $A \in \mathscr{B}(H).$ \\
Conversely, suppose for all $A \in \mathscr{B}(H),$ $\Phi(A) = \lambda UAU^*,$ where $U$ is unitary and $|\lambda|=1.$ Then $$(\Phi(A)-\Phi(B))\Phi(C) = \lambda^2 U(A-B)CU^*.$$ Therefore,
$q((\Phi(A)-\Phi(B))\Phi(C)) = q((A-B)C)$ for all $A,B, C \in \mathscr{B}(H).$
\end{proof}
\begin{theorem}\label{thm 3.10}
Let $H$ be an infinite dimensional complex Hilbert space and $\Phi: \mathscr{B}(H) \to \mathscr{B}(H)$ be surjective linear map. Then 
$$m((\Phi(A)- \Phi(B))\Phi(C)) = m((A-B)C),$$ for all $A,B,C \in \mathscr{B}(H)$ if and only if there exists a unitary operator $U \in \mathscr{B}(H)$ and a complex scalar $\lambda$ with $|\lambda| = 1$ such that $\Phi(A) = \lambda UAU^*$ for all $A \in \mathscr{B}(H).$
\end{theorem}
\begin{proof}
Taking $B = -A,$ we get $m(\Phi(A)\Phi(C)) = m (AC)$ for all $A, C \in \mathscr{B}(H).$ Recall the fundamental identity $q(T) = m(T^*)$ for all $T \in \mathscr{B}(H).$ Define $\Psi: \mathscr{B}(H) \to \mathscr{B}(H)$ by $\Psi(A) = (\Phi(A^*))^*$ for all $A \in \mathscr{B}(H).$ It is evident that $\Psi$ is linear. Since $\Phi $ is surjective. So $\Psi$ is also surjective. Moreover, $\Psi(A)\Psi(C) = (\Phi(C^*)\Phi(A^*))^*.$ Consequently, $q(\Psi(A)\Psi(C)) = m(\Phi(C^*)\Phi(A^*)) = q(AC)$ for all $A, C \in \mathscr{B}(H).$ Using Theorem \ref{thm 3.8}, there exists a unitary operator $U \in \mathscr{B}(H)$ and a complex scalar $\lambda$ with $|\lambda|= 1$ such that $\Psi(A) = \lambda UAU^*$ for all $A \in \mathscr{B}(H).$ Therefore, $\Phi(A) = \overline{\lambda}UAU^*$ for all $A \in \mathscr{B}(H).$\\
Conversely, Suppose $\Phi(A) = \lambda UAU^*,$ for all $A \in \mathscr{B}(H)$ and $|\lambda|=1$ with unitary $U \in \mathscr{B}(H).$ Then $$m((\Phi(A)-\Phi(B))\Phi(C)) = |\lambda|^2m((A-B)C) = m((A-B)C),$$ for all $A, B, C \in \mathscr{B}(H).$
\end{proof}

In \cite{brevsar1996linear}, the authors characterized surjective linear maps on \(\mathscr B(\mathscr X),\) where $\mathscr{X}$ is a complex Banach space,  that preserve the spectral radius, showing that every such map is a scalar multiple of a similarity or an adjoint-type transformation. Recently, Aznay, and Chabbabi \cite{azna2026surjective} established a structural characterization of surjective maps $\Phi$ on $\mathscr{B(\mathscr X)},$ where $\mathscr X$ is a infinite dimensional complex Banach space, preserving the spectral radius of operator products, namely, 
$$r(\Phi(A)\Phi(B)) = r(AB),$$ for all $A, B \in \mathscr{B}(\mathscr X).$ Moreover, In \cite{hou2007linear}, the authors characterized linear maps on \(\mathscr B(H),\) where $H$ is infinite-dimensional, preserving essential spectral properties, Fredholm and semi-Fredholm operators, and the closeness of operator ranges. Theorem \ref{thm 3.14} provides a structural characterization of surjective linear maps  $\Phi: \mathscr{B}(H) \to \mathscr{B}(H),$ where $H$ is infinite dimensional, that preserve the spectral radius of operator products involving differences, namely,
$$r((\Phi(A)-\Phi(B))\Phi(C)) = r((A-B)C),$$ for all $A,B, C \in \mathscr{B}(H).$ Before that, we establish a lemma needed for the proof of Theorem \ref{thm 3.14}.
\begin{lemma}\label{lemma 3.13}
Let $A, B \in \mathscr{B}(H).$ Suppose that $r(AC) = r(BC)$ for all $C \in \mathscr{B}(H).$ Then either $A=B=0$ or there exists a non zero $\lambda \in \mathbb{C}$ with $|\lambda|=1$ such that $A = \lambda B.$
\end{lemma}
\begin{proof}
For all $x, y \in H,$ we have 
\begin{equation}\label{equ 4}
r(A(x \otimes y)) = |\langle Ax, y\rangle| = r(B(x \otimes y)) = |\langle Bx, y\rangle|.
\end{equation}
Now we consider two cases.\\
$\textit{Case(I)}:$ For all $x, y \in H,$ $|\langle Ax, y\rangle| = 0$ if and only if $|\langle Bx, y\rangle| =0.$ Consequently, $A=B.$\\
$\textit{Case(II)}:$ Suppose $B\ne 0.$ For each $Bx \ne 0,$ we  get $Ax = \mu_x Bx$ with $|\mu_x|=1$ because $\|Ax\|\|Bx\| \ge |\langle Ax, Bx\rangle| = \|Bx\|^2$ and $\|Bx\| \|Ax\| \ge |\langle Bx, Ax\rangle| = \|Ax\|^2.$  Using equation (\ref{equ 4}), when $Bx =0$ then $Ax =0.$\\
We first consider the case in which $\dim R(B)\ge 2.$ Choose $x_1, x_2 \in H$ such that $Bx_1, Bx_2$ are linearly independent. By linearity, we obtain 
$$\mu_{x_1 + x_2}B(x_1 + x_2) = A(x_1 +x_2) = \mu_{x_1}Bx_1 + \mu_{x_2}Bx_2.$$ Since $Bx_{1}, Bx_{2}$ are linearly independent. Then $\mu_{x_1 + x_2} = \mu_{x_1} = \mu_{x_2}.$ By repeating this argument, and connecting any two non zero vectors in $R(B)$ through linearly independent vectors when necessary, we conclude that $\mu_x$ is independent of $x \in H.$Thus there exists $\lambda \in \mathbb{C}$ with $|\lambda| =1$ such that $A = \lambda B.$\\
Finally, suppose that $\dim R(B) =1.$ Since $Ax$ and $Bx$ are linearly dependent for all $x \in H.$ Then we have two linear functional $f, g$ from $H$ to $\mathbb{C}$ such that $$Ax = f(x)u \text{ and } Bx = g(x) u,$$ for some non zero vector $u \in R(B)$ with $|f(x)| = |g(x)|$ for all $x \in H.$ There exists a nonzero $x_0 \in H$ such that $g(x_0) \ne 0.$  So $
|f(x_0)|=|g(x_0)|\ne0.$ Define
\[
\mu=\frac{f(x_0)}{g(x_0)}.
\]

Clearly, \(|\mu|=1\). Now, for arbitrary \(x\in H\) and \(t\in\mathbb C\), we obtain
$$
|f(x_0)+tf(x)|^2
=
|g(x_0)+tg(x)|^2.
$$ Expanding both sides gives
\[
\begin{aligned}
&|f(x_0)|^2
+2\operatorname{Re}\big(t f(x)\overline{f(x_0)}\big)
+|t|^2|f(x)|^2\\
&=
|g(x_0)|^2
+2\operatorname{Re}\big(t g(x)\overline{g(x_0)}\big)
+|t|^2|g(x)|^2.
\end{aligned}
\]

Using the property, $|f(x)|=|g(x)|$ for all $x \in H,$ we get
$$\operatorname{Re}\left(
t\left[f(x)\overline{f(x_0)}
-g(x)\overline{g(x_0)}\right]\right)=0,$$ for every \(t\in\mathbb C\). Taking \(t=1\) and \(t=i\), we conclude that
$f(x)\overline{f(x_0)}
=g(x)\overline{g(x_0)}.$ Since \(f(x_0)=\mu g(x_0)\) and \(|\mu|=1\), it follows that $f(x)=\mu g(x).$
As \(x\in H\) was arbitrary, we obtain ${f=\mu g, ~ |\mu|=1.}$ Therefore, $A = \mu B$ for some unit element $\mu \in \mathbb{C}.$

\end{proof}
\begin{theorem}\label{thm 3.14}
Let $H$ be an infinite dimensional complex Hilbert space and $\Phi: \mathscr{B}(H) \to \mathscr{B}(H)$ be a surjective  map. Then 
$$r((\Phi(A)- \Phi(B))\Phi(C)) = r((A-B)C),$$ for all $A,B,C \in \mathscr{B}(H)$ if and only if there exists an invertible linear or conjugate linear map  $T$ on $H$ and a complex scalar $\lambda \in \mathbb{C}$ with $|\lambda| =1$ such that $\Phi$ has one of the following forms: 
\begin{enumerate}
\item $\Phi(A) = \lambda TAT^{-1}$ for all $A \in \mathscr{B}(H),$
\item $\Phi(A) = \lambda TA^*T^{-1}$ for all $A \in \mathscr{B}(H).$
\end{enumerate}
\end{theorem}
\begin{proof}
By surjectivity, there exists $C_0 \in \mathscr{B}(H)$ such that $\Phi(C_0) =0.$ Taking $C=C_0$ in the main equation, we get $r((A-B)C_0) =0$ for all $A, B \in \mathscr{B}(H).$ Now consider, $A-B = C_0^*,$ we obtain $0= r(C_0^*C_0) = \|C_0^*C_0\| = \|C_0\|^2.$ So $\Phi(0) =0.$ Next we show that $\Phi$ is bijective.\\
Suppose that $\Phi(A) = \Phi(B).$ It follows that $r((A-B)C) = 0$ for all $C \in \mathscr{B}(H).$ Choosing $C= (A-B)^*,$ we obtain $\|A-B\|^2 =0,$ hence $A=B.$ Since $\Phi$ is surjective. Consequently, $\Phi$ is bijective. Now putting $B=0$ in the main equation and applying $\Phi(0) =0,$ we have \begin{equation}\label{equ 5}
r(\Phi(A)\Phi(C)) = r(AC),
\end{equation}
for all $A, C \in \mathscr{B}(H).$ It follows that 
\begin{equation}\label{equ 6}
r(\Phi(A-B)\Phi(C)) = r((A-B)C) = r((\Phi(A)-\Phi(B))\Phi(C)).
\end{equation}
Since $\Phi$ is surjective. Applying Lemma \ref{lemma 3.13}, we obtain $$\Phi(A)-\Phi(B) = \mu(A,B) \Phi(A-B),$$ for all $A,B \in \mathscr{B}(H)$ and $|\mu(A,B)|=1$ whenever $A \ne B.$\\
Now, we invoke Theorem 1.1 \cite{azna2026surjective}, which gives precisely the following alternatives:
\begin{equation}\label{equ 7}
\Phi(A) = \lambda(A)SAS^{-1}, \text{ for all } A \in \mathscr{B}(H),
\end{equation}
\begin{equation}\label{equ 8}
\Phi(A) = \lambda(A)SA^*S^{-1}, \text{ for all } A \in \mathscr{B}(H),
\end{equation}
where $\lambda: \mathscr{B}(H) \to \mathbb{C}\setminus \{0\}$ and $S$ a is bounded invertible, linear or conjugate linear map on $H.$\\
Now first consider the representation (\ref{equ 7}). Since $$r(A^*C^*) = r(CA)^* = r(CA) = r(AC)$$ for all $A, C \in \mathscr{B}(H).$ Using the identity (\ref{equ 5}), we have 
$$r(AC) = r(\Phi(A)\Phi(C)) = |\lambda(A)| |\lambda(C)|R(SACS^{-1}) = |\lambda(A)| |\lambda(C)| r(AC),$$ for all $A, C \in \mathscr{B}(H).$ So $|\lambda(A)| |\lambda(C)| =1$ whenever $r(AC) > 0.$ Taking $A= I=C,$ we obtain $|\lambda(I)|= 1. $\\
Now let $A \ne 0.$ Choose $x \in H$ such that $Ax \ne 0.$ By choosing $y \in H$ such that $\langle x, y \rangle \ne 0,$ $\langle Ax, y \rangle \ne 0$ and setting $C= x \otimes y,$ we have $$r(IC) =r(C) = |\langle x, y \rangle| >0, \quad \text{and} \quad r(AC) = |\langle Ax,y\rangle| >0.$$ Now, it is evident that $|\lambda(C)| = |\lambda(I)||\lambda(C)| =1,$ hence  $1= |\lambda(A)||\lambda(C)| = |\lambda(A)|.$ Moreover, 
$$|\lambda(A)| =1, \quad \text{ for all} ~A\in \mathscr{B}(H)\setminus\{0\}.$$ Now we show that the function $\lambda$ is constant. From equation (\ref{equ 7}) and the relation $\Phi(A)-\Phi(B) = \mu(A,B) \Phi(A-B),$ for all $A,B \in \mathscr{B}(H)$ and $|\mu(A,B)|=1$ whenever $A \ne B,$ we get 
$$\lambda(A)A - \lambda(B)B = \mu(A, B) \lambda(A-B) (A-B),$$ for all $A, B \in \mathscr{B}(H).$\\
Let $A, B \in \mathscr{B}(H)$ be linearly independent. Comparing the coefficient, we obtain 
$\lambda(A) = \mu(A,B) \lambda(A-B) = \lambda(B).$ Since $\mathscr{B}(H)$ has dimension greater than 1, any two non zero operators can be compared through a third operator linearly independent with both. It follows that the function $\lambda$ is constant, say $\gamma,$ on $\mathscr{B}(H)$ with $|\gamma|=1.$ The same argument applies to equation (\ref{equ 8}).  Therefore, there exists a fixed $\gamma \in \mathbb{C}$ with $|\gamma|=1$ and a bounded invertible linear or conjugate linear $S$ map on $H$ such that $\Phi$ has one of the following forms: 
\begin{enumerate}
\item $\Phi(A) = \gamma SAS^{-1}$ for all $A \in \mathscr{B}(H),$
\item $\Phi(A) = \gamma SA^*S^{-1}$ for all $A \in \mathscr{B}(H).$
\end{enumerate}
Conversely, suppose that $\Phi(A) = \lambda TAT^{-1}$ for all $A \in \mathscr{B}(H),$ where $T$ is bounded invertible, linear or conjugate linear with $|\lambda| =1.$ Then it is obvious that $r((\Phi(A)- \Phi(B))\Phi(C)) = r((A-B)C),$ for all $A,B,C \in \mathscr{B}(H).$\\
Suppose that $\Phi(A) = \lambda TA^*T^{-1}$ for all $A \in \mathscr{B}(H),$ where $T$ is bounded invertible, linear or conjugate linear with $|\lambda| =1.$ Then $$r((\Phi(A)- \Phi(B))\Phi(C)) = r((A-B)^*C^*) = r(C(A-B)) = r((A-B)C),$$ for all $A, B, C \in \mathscr{B}(H).$

\end{proof}

In \cite{bourhim2025unitary}, the authors characterized all bijective linear maps \(\Phi:\mathscr{B}(H)\to\mathscr{B}(H)\) satisfying
\[
w(\Phi(A))=w(\Phi(B)),\quad \text{ for all }A,B\in\mathscr{B}(H),
\]

when \(H\) is an infinite-dimensional separable Hilbert space. Recently, Amara \cite{amara2026spectral} extended this characterization to both finite- and infinite-dimensional Hilbert spaces, without requiring separability. In the next theorem, we characterize all bijective linear maps \(\Phi:\mathscr{B}(H)\to\mathscr{B}(H)\) satisfying $\Phi(I) =I$ and   \[
w_0(\Phi(A))=w_0(\Phi(B)),\qquad A,B\in\mathscr{B}(H),
\]
 when \(H\) is an infinite-dimensional Hilbert space, not necessarily separable. We first establish a lemma needed for the proof of Theorem \ref{thm 3.12}.\\
 Define a function $h: \mathscr{B}(H) \to \mathbb{R}$ by 
 $$h(A) = \sup_{\|x\| =1} \Re\langle Ax, x\rangle,$$ for all $A \in \mathscr{B}(H).$

\begin{lemma}\label{lemma 3.11}
$h(A) = \lim_{r \to \infty} (w_0(rI+A)-r) \text{ for all } A \in \mathscr{B}(H).$
\end{lemma}
\begin{proof}
Fix $r > \|A\|.$ By the definition of the maximal numerical range, each $z \in W_0(rI +A)$ can be written as $z = r +\lambda,$ where there is a sequence of unit vectors $\{x_n\}$ satisfying $$\|(rI + A)x_n\| \to \|rI +A\| \quad \text{and} \quad \langle Ax_n, x_n\rangle  \to \lambda$$ as $n \to \infty.$ Using $\|(rI+ A)x_n\|^2 = r^2 + 2r\Re\langle Ax_n, x_n\rangle + \|Ax_n\|^2,$ we obtain, on taking limits in both side
$$\|rI+A\|^2 = r^2 +2r \Re \lambda + l,$$ where $0 \le l = \lim_{n \to \infty}\|Ax_n\|^2 \le \|A\|^2,$ if necessary take subsequence of $\{x_n\}$ then renaming as $\{x_n\}.$ On the other hand, by the definition of $h(A),$ we get
\begin{align*}
r^2 + \|A\|^2 + 2r \Re \lambda\ge r^2 + 2r \Re \lambda + l &= \|rI +A\|^2\\ &= \sup_{\|x\|=1} r^2 + 2r\Re\langle Ax, x\rangle + \|Ax\|^2\\&\ge r^2 + 2rh(A).
\end{align*}
Consequently, 
\begin{equation}\label{equ 2}
h(A)-\frac{\|A\|^2}{2r} \le \Re \lambda \le h(A).
\end{equation} 
Since $|\lambda| \le \|A\|$ and $r > \|A\|.$ Then $r + \Re \lambda >0.$ Hence
\begin{align*}
r + \Re \lambda \le |r + \lambda|\le \sqrt{(r +\Re \lambda)^2 + (\operatorname{Im} \lambda)^2} \le (r + \Re \lambda) + \frac{(\operatorname{Im} \lambda)^2}{2(r + \Re \lambda)}.
\end{align*}
By using equation (\ref{equ 2}), we obtain 
$$r + h(A)-\frac{\|A\|^2}{2r} \le |r+\lambda| \le r+h(A) + \frac{\|A\|^2}{2(r -\|A\|)}.$$
Since $z \in W_0(rI +A),$ taking the supremum over $|z|,$ we get $$r + h(A)-\frac{\|A\|^2}{2r} \le w_0(rI+A) \le r+h(A) + \frac{\|A\|^2}{2(r -\|A\|)}.$$ Letting $r \to \infty,$ we get the required formula $$h(A) = \lim_{r \to \infty} (w_0(rI+A)-r).$$
 \end{proof}
\begin{theorem}\label{thm 3.12}
Let $H$ be infinite-dimensional and $\Phi: \mathscr{B}(H) \to \mathscr{B}(H)$ be bijective linear map with $\Phi(I) =I$. 
The following statements are equivalent:
\begin{enumerate}
\item \(w_0(\Phi(T))=w_0(\Phi(S))\) whenever \(T\sim_u S\).
\item There exists a unitary operator \(U\in\mathscr{B}(H)\) such that \(\Phi\) has one of the following forms:
\[
\Phi(A)= UAU^*
    \quad\text{or}\quad
    \Phi(A)= UA^{\mathrm{tr}}U^*, \text{ for all}~ A\in\mathscr{B}(H).
\]
\end{enumerate}
\end{theorem}
\begin{proof} $(1) \implies (2)$
For a unitary operator $V \in \mathscr{B}(H),$ we have $$VTV^*+ \lambda I = V(T+\lambda I)V^*,$$ for all $T \in \mathscr{B}(H).$ From our assumptions give $$w_0(\Phi(VTV^*)+ \lambda I) = w_0(\Phi(VTV^* + \lambda I)) = w_0(\Phi(T+\lambda I)) = w_0(\Phi(T)+ \lambda I),$$ for every $\lambda \in \mathbb{C}.$ Consider $A = \Phi(T), \quad B=\Phi(VTV^*), \quad \text{and} \quad \lambda =r.$ Using Lemma \ref{lemma 3.11} and limiting $r \to \infty$ give $h(\Phi(T)) = h(A) =h(B)= h(\Phi(VTV^*)).$ Similarly, we can show that $h(e^{-i\theta}A) = h(e^{-i\theta}B)$ for all $\theta \in [0, 2\pi).$ But $h(e^{-i\theta}A) = \sup_{\|x\|=1} \Re (e^{-i\theta}\langle Ax,x\rangle).$ So the support function of $W(A)$ is direction $\theta$ (see \cite{bourdon2000numerical}) is $$S_{W(A)}(\theta) = \sup_{z \in W(A)}\Re(e^{-i\theta}z) =h(e^{-i\theta}A) = h({e^{-i\theta}B}) = S_{W(B)}(\theta),$$ for all $\theta \in [0,2\pi).$ Consequently, $$\overline{W(A)} = \bigcap_{\theta \in [0, 2\pi)}\big\{z: \Re(e^{-i\theta}z) \le S_{W(A)}(\theta)\big\} = \overline{W(B)}.$$ Hence $\overline{W(\Phi(T))} = \overline{W(\Phi(VTV^*))}.$ It follows that $w(\Phi(T)) = w(\Phi(VTV^*)).$ Using Theorem 2.4 \cite{amara2026spectral}, there exists a unitary operator $U \in \mathscr{B}(H)$ and non zero scalar $c \in \mathbb{C}$ such that either $\Phi(A) = cUAU^*$ for all $A \in \mathscr{B}(H)$ or $\Phi(A) = cUA^{\mathrm{tr}}U^*$  for all $A \in \mathscr{B}(H).$  Since $\Phi(I) =I.$ Therefore, There exists a unitary operator \(U\in\mathscr{B}(H)\) such that \(\Phi\) has one of the following forms:
\[
\Phi(A)= UAU^*
    \quad\text{or}\quad
    \Phi(A)= UA^{\mathrm{tr}}U^*, \text{ for all}~ A\in\mathscr{B}(H).
\]
$(2) \implies (1)$ It is obvious by using the property $$w_0(T^{\mathrm{tr}}) = w_{0}(T) = w_{0}(UTU^*),$$ for all $T \in \mathscr{B}(H)$ and $U \in \mathscr{B}(H)$ is unitary.
\end{proof}

\section*{Funding}
The author gratefully acknowledges support from GA~\v{C}R(Grant No.~25-15444K) and RVO:67985840. 

\section*{Data Availability}
Data sharing is not applicable to this article as no data sets were generated or analyzed during the current study.

\section*{Declarations}
  The author confirms that there are no conflicts of interest regarding the publication of this manuscript.

\end{document}